\documentclass{amsart} 
\usepackage{epsf} 
\usepackage{amssymb,latexsym, amsmath, amscd, array, graphicx} 
\def\hepsffile{\leavevmode\epsffile} 

\swapnumbers 
\numberwithin{equation}{section} 

\theoremstyle{plain} 
\newtheorem{thm}{Theorem}[section] 
\newtheorem{cor}[thm]{Corollary} 
\newtheorem{lem}[thm]{Lemma} 
 
\newtheorem{prop}[thm]{Proposition} 

\newcommand\theoref{Theorem~\ref} 
 
\newcommand\propref{Proposition~\ref}

\theoremstyle{definition} 
\newtheorem{defin}[thm]{Definition}

\newtheorem{rem}[thm]{Remark} 
\newtheorem{ex}[thm]{Example} 
\newtheorem{con}[thm]{Construction} 
 
\newtheorem{sourex}[thm]{Source of Examples of 
propagations without dangerous intersections}

\def\Im{\protect\operatorname{Im}}

\def\Ker{\protect\operatorname{Ker}}

\def\sec{\protect\operatorname{sec}} 
\def\sign{\protect\operatorname{sign}}

\def\lk{\protect\operatorname{AL}} 
 
\def\pr{\protect\operatorname{pr}} 
\def\cri{\protect\operatorname{CR}} 
\def\aff{\protect\operatorname{aff}} 
\def\clk{\protect\operatorname{lk}} 
\def\exp{\protect\operatorname{exp}} 
\def\det{\protect\operatorname{det}} 

\def\eps{\varepsilon} 

\def\Z{{\mathbb Z}} 
\def\Q{{\mathbb Q}} 
\def\R{{\mathbb R}} 
\def\N{{\mathbb N}} 
\def\A{{\pmb A}}

\def\e{{\mathbf e}} 

\def\1{\hbox{\rm\rlap {1}\hskip.03in{\rom I}}} 
\def\Bbbone{{\rm1\mathchoice{\kern-0.25em}{\kern-0.25em} 
{\kern-0.2em}{\kern-0.2em}I}} 

\def\pp{\medskip{\parindent 0pt \it Proof.\ }} 
\def\wt{\widetilde} 
\def\wh{\widehat} 
\def\m{\medskip} 
\def\ov{\overline} 

\begin{document} 
\date{January 23, 2004} 
\leftline{ } 
\centerline{ } 
\title[Affine Linking Numbers and Causality Relations for 
Wave Fronts] 
{Affine Linking Numbers and Causality Relations for Wave 
Fronts} 
\author[V.~Chernov (Tchernov) and Yu.~Rudyak]{Vladimir V. 
Chernov 
(Tchernov) and Yuli B. Rudyak} 
\address{V. Chernov, Department of Mathematics, 
6188 Bradley Hall, Dartmouth College, Hanover NH 03755, 
USA} 
\email{Vladimir.Chernov@dartmouth.edu} 
\address{Yu. Rudyak, Department of Mathematics, University 
of Florida, 358 
Little Hall, Gainesville, FL 32611-8105, USA} 
\email{rudyak@math.ufl.edu} 



\Large 

\begin{abstract} 
Two wave fronts $W_1$ and $W_2$ that 
originated at some points of the manifold $M^n$ are said 
to be causally related 
if one of them passed through the origin of the other 
before the other appeared. 
We define the causality 
relation invariant $\cri(W_1, W_2)$ to be the algebraic 
number 
of times the earlier born front passed through the origin 
of the other front 
before the other front appeared. Clearly, if $\cri(W_1, 
W_2)\neq 0$, then $W_1$ 
and $W_2$ are 
causally related. If $\cri(W_1, W_2)= 0$, then we 
generally can not make any 
conclusion about fronts being causally related. However we 
show that for front 
propagation given by a complete Riemannian metric of non- 
positive sectional 
curvature, $\cri(W_1, W_2)\neq 0$ if and only if the two 
fronts are causally 
related. The models where the law of propagation is given 
by a metric of 
constant sectional curvature are the famous Friedmann 
Cosmology models. 

The classical linking number $\clk$ is a $\Z$-valued 
invariant of two zero 
homologous submanifolds. We construct the affine linking 
number generalization 
$\lk$ of the $\clk$ invariant to the case of linked $(n-1) 
$-spheres in the total 
space of the unit sphere tangent bundle $(STM)^{2n-1}\to 
M^n$. 
For all $M$, except of odd-dimensional rational homology 
spheres, $\lk$ allows 
one 
to calculate the value of $\cri(W_1, W_2)$ from the 
picture of the two wave 
fronts at a certain moment. This calculation is done 
without the knowledge of 
the front propagation law and of their points and times of 
birth. 
Moreover, in fact we even do not need to know the 
topology of $M$ outside of a part $\overline M$ of $M$ 
such that 
$W_1$ and $W_2$ are null-homotopic in $\overline M$. 
\end{abstract}

\maketitle

\section*{Introduction} 
In this paper the word ``smooth'' means $C^{\infty}$. 
Throughout this paper $M$ is a smooth connected oriented 
Riemannian $n$-dimensional manifold (not 
necessarily compact). 
In our paper the wave fronts on $M$ are assumed to be 
parametrized by smooth 
mappings of manifolds. They are not assumed to be immersed 
and are allowed to 
have various singularities including degenerations of the 
differential. The popular case of wave fronts being 
projections to $M$ of 
Legendrian mappings of $N^{n-1}$ (parameterizing the 
fronts) to the unit 
cotangent bundle $(ST^*M)^{2n-1}\subset T^*M=TM$ of $M$ is 
easily obtained from 
our paper 
as a particular, though a very important case. (The front 
stays smooth, though 
not immersed, at the cusp points and other singularities 
of the projections to 
$M$ of Legendrian mappings $N^{m-1}\to (ST^*M)^{2n-1}$.)

Let $W_1$ and $W_2$ be two wave fronts which are 
propagating in $M$. 
(Generally, we assume that the fronts have different 
propagation laws.) 
We define a {\it dangerous intersection} between 
the fronts $W_1(t)$ and $W_2(t)$ at some moment of time 
$t$ to be a point $x$ 
where the fronts intersect and have the same direction of 
propagation (see Section \ref{prelim} for the precise 
definition). 

A passage of the front $W_1(t)$ through the birth point 
of the front $W_2$ at the moment of time $t$ 
before the front $W_2$ originated is called the {\it baby- 
intersection}.

It turns out that we can associate to each dangerous 
intersection as well as 
to each baby-intersection a sign (i.e. a number $\pm 1$ ). 
The 
sum of the signs up to 
a moment $t$ is called a {\it causality relation 
invariant} 
and is denoted by 
$\cri(W_1(t),W_2(t))$. 

In particular, assume that there are no dangerous 
intersections 
of $W_1(t)$ and $W_2(t)$ for all $t$. (Many such examples 
are constructed in 
Section \ref{prelim} and, for example, light and sound 
wave fronts have this 
property.) Then the causality relation invariant tells 
us the algebraic number of times the earlier-born wave 
front passed through the 
birth point 
of the other front before the other front originated. 

Two wave fronts $W_1$ and $W_2$ that 
originated at some points of the manifold $M^m$ are said 
to be {\em causally 
related\/} 
if one of them passed through the origin of the other 
before the other appeared. 
Clearly, if $\cri(W_1, W_2)\neq 0$ and no dangerous 
intersections occurred 
during the propagation, 
then $W_1$ and $W_2$ are causally related. 

If $\cri(W_1, W_2)= 0$, then we generally can not make any 
conclusion about 
fronts being causally related. However we show that in 
case of front propagation 
given by a complete Riemannian metric of non-positive 
sectional curvature, 
$\cri(W_1, W_2)\neq 0$ if and only if the two fronts are 
causally related, 
see~\ref{Friedmann}. The models where the law of 
propagation is given by a 
metric of constant sectional curvature are the famous 
Friedmann Cosmology 
models.

We are interested in reconstructing the value $\cri(W_1 
(t),W_2(t))$ from the current shape of the wave fronts 
only, 
without the knowledge of the propagation laws, of the 
birth-points 
of the fronts, of the topology of $M$ etc.

It turns out that, having the current picture only, we can 
evaluate 
$\cri(W_1(t),W_2(t))\in \Z$ modulo a certain $m\in \Z$ 
that depends on $M$. This $m$ 
is zero if $M$ is not an odd-dimensional rational homology 
sphere, 
see~\propref{usefulprop}, and when $m=0$ we can completely 
reconstruct the value 
of $\cri(W_1(t),W_2(t))\in \Z$. Furthermore, for $M$ odd-dimensional 
this $m$ is 
divisible by the order of $\pi_1(M)$, see 
\propref{usefulprop}. In 
particular, if 
$\pi_1(M)$ is infinite then $m=0$, i.e. we can completely 
evaluate $\cri$ from 
the current picture. The really bad case $m=1$ (when 
we can not say anything about $\cri$) appears only when 
$M$ is an 
odd-dimensional homotopy sphere. In particular, for 
Friedmann cosmology models 
based on a 
Riemannian metric of non-positive sectional curvature our 
methods allow one to 
detect precisely whether two wave fronts are causally 
related or not from the 
current picture of the wave fronts only,~see~\ref 
{Friedmannlinking}. 

To evaluate $\cri$ modulo $m$, we introduce an invariant 
$\lk 
\in \Z/m\Z$, the so-called affine linking invariant, which 
depends on the current picture only. Then we notice that 
$\cri$ and $\lk$ are 
congruent modulo $m$. Here the biggest technical 
difficulty 
appears, since in 
order to define $\lk$ we must define the ``linking 
number'' 
for two spheres that are non-homologous to zero. 
(This affine linking number can be shown to be a 
particular case of a very general affine linking invariant 
discussed in our later work~\cite 
{ChernovRudyakGeneralLinking}. However, results of this 
paper are independent from the results 
of~\cite{ChernovRudyakGeneralLinking} and are not 
corollaries of our results obtained later.)

This theory has the following physical interpretation. Let 
$\overline M$ be the part 
of the manifold (the universe) $M$ such that $\overline M$ 
contains the current 
picture of wave fronts $W_1, W_2$, and $W_1, W_2$ are 
contractible in $\overline 
M$. 

We transform the wave fronts via certain allowable moves 
to 
trivial fronts, i.e. 
small spherical fronts with the canonical orientation and 
coorientation, located 
far away from each other. The allowable moves should be 
thought of as 
generalized Reidemeister moves: they are the passages 
through generic 
singularities (in both directions) of wave fronts and 
dangerous intersection moves. 
We count the change of the invariant $\cri$ that occurs in 
the process of this 
formal deformation, and it turns out that this change is 
congruent modulo $m$ 
with the (unknown!) value $\cri(W_1,W_2)$ of the current 
picture. In particular, 
as we have already mentioned, if $M$ is not an odd- 
dimensional rational homology 
sphere or if $\pi_1(M) 
$ is infinite, then 
we can completely compute $\cri$ from the current picture, 
without any knowledge of the propagations, moments and 
points of birth 
of the fronts, and topology of $M$ 
outside of $\overline M$.

The following observation seems to be interesting. 
Suppose that we have two pictures 
of two pairs of fronts $(W_1, W_2)$ and $(W'_1, 
W'_2)$ made at 
two unknown moments of time $t_0$ and $t_1$. (We assume 
that 
both pairs are free of dangerous intersection points.) 
Assume 
that we know that the propagation laws for the two fronts 
are such that 
the dangerous intersection points cannot appear during the 
propagation 
and that 
$\cri(W_1(t_0), W_2(t_0))$ and $\cri (W'_1(t_1), W'_2 
(t_1)) 
$ are not comparable modulo $m$. 
Then we can conclude that the pairs $(W_1(t_0), W_2(t_0))$ 
and $(W'_1(t_1), W'_2(t_1))$ of wave fronts are 
not the pictures of the same pair of fronts taken at 
different 
moments of time.

Note that in these calculations we disregard the 
dangerous self-intersections of wave fronts. (In a sense 
this is similar to 
the theory of link homotopy where different components of 
links are not 
allowed to intersect through possible deformations, but 
self-intersections 
are allowed.) The study of self-intersections of fronts on 
surfaces was 
initiated by the ground breaking work of Arnold~\cite 
{Arnold}, see also 
\cite{Aicardi1, Aicardi2, ChmutovGoryunov, Goryunov1, 
Goryunov2, Inshakovcurves, 
Inshakovgroupscurves, Polyak, 
PolyakBennequin, Chernovfronts, Chernovcurves, 
Chernovamsvolume, Tabachnikov}. 
The methods developed in this paper allow us to calculate 
the 
algebraic number of dangerous self-intersection points 
that arise under the 
propagation of fronts on manifolds of arbitrary 
dimensions, we do 
it in a next paper.

The following physical speculations related to the $\cri$ 
invariant seem to be possible. Assume that the space-time 
is 
topologically a product $M^n\times\R$, and that the 
observable 
universe $\overline M$ is so big that we are not able to 
see the 
current picture of wave fronts (due to the finiteness of 
the speed 
of light). The propagating fronts define the mapping of 
the 
cones $C_1, C_2 $ (over the sphere $S^{n-1}$ 
parameterizing the 
fronts at every moment of time) into $M \times \R$. Let 
$\sec: 
\overline M \to \overline M\times \R$ be a section of the 
projection $p_ {\overline M}: \overline M\times \R \to 
\overline 
M$, and let $\overline W_i=C_i \cap \sec (\overline M), 
i=1,2$. Assume that 
the law of propagation is such that dangerous 
intersections do not 
occur. Then, similarly to the above, we can restore the 
number of 
baby-intersections from the picture of images of 
$\overline W_i$ under the 
projection $p_ {\overline M}$. The section $\sec$ can be 
thought of 
as the picture of the universe that we see as the light 
from the 
points of $\overline M$ reaches the observer, and thus 
$\overline W_i$ can 
be regarded as the picture of fronts that we actually see.

Low \cite{Low1, Low2, Low3} attacked a similar problem for 
$M=\R^3$, where he considered a linking invariant for 
linked 
cones $C_1, C_2$ as above. In this case the linking 
numbers can be 
constructed directly via the approach of 
Tabachnikov~\cite{Tabachnikov}, because $\R^3$ has the 
topological 
end. (Some very interesting results relating causality for 
fronts on $\R^2$ and linking were 
obtained recently by Natario and Todd~\cite{NatarioTod}.)

The paper is organized as follows. In Section~\ref{prelim} 
we 
discuss some preliminary information, in Section~\ref{cri} 
we 
define the invariant $\cri$, in Section~\ref{hom} we prove 
homotopy theoretical results which we use in order to 
define the 
invariant $\lk$, in Section~\ref{lk} we define the 
invariant $\lk$ 
and state the relation between $\cri$ and $\lk$, in 
Section~\ref{mod} 
we treat the case of propagation with respect to a certain 
Riemannian 
metric, in Section~\ref{examples} we give some examples 
and applications.

\section{Preliminaries: propagation laws, propagations and 
dangerous intersections}\label{prelim} 

We denote by $p_T: TM \to M$ the tangent 
bundle over $M$. Let $s: M\to TM$ be the zero section of 
the tangent bundle. We set $\overline {TM}=TM\setminus s(M) 
$. 
The multiplicative group $\R^+$ of positive real 
numbers acts fiberwise on $TM\setminus s(M)$ by 
multiplication, and we set 
\begin{equation*} 
STM=(TM\setminus s(M))/\R^+. 
\end{equation*} 

Let $\overline p: \overline {TM} \to STM$ be the quotient 
map. Clearly, 
the projection $p_T: TM \to M$ yields the commutative 
diagram 
$$ 
\CD 
TM @>\supset>> \overline {TM} @>\overline p>> STM\\ 
@Vp_TVV @VVV @VV\pr V\\ 
M @= M @= M 
\endCD 
$$ 
It is easy to see that $\pr: STM \to M$ is a locally 
trivial 
bundle with the fiber 
$S^{n-1}$, we call this bundle the {\it spherical tangent 
bundle}. 

Given $x\in M^n$, we denote by $S^{n-1}_x$ (or just by 
$S_x$) the 
fiber $\pr^{-1}(x)$ over $x$ of the spherical tangent 
bundle and 
by $T_xM$ the tangent space to $M$ at $x$. 

Since $M$ is orientable, the bundle $\pr: STM \to M$ is 
also 
orientable, and in order to orient $STM$ it suffices to 
orient the 
fiber $S^{n-1}_x$. We do it as follows. Choose an 
orientation 
preserving chart for $M$ centered at $x$ and let $S$ be a 
small 
$(n-1)$-sphere centered at $x$. We equip $S$ with the 
unique 
orientation $o$ by requiring that the pair ($o$, outer 
normal 
vector to $S$) gives us the orientation of $M$. 

Given $s\in S$, the radius-vector from $x$ to $s$ can be 
regarded 
as a nonzero tangent vector to $M$ at $x$, i.e., as a 
point of 
$S^{n-1}_x$. In this way we get a diffeomorphism $\psi: 
S\to S_x$ 
which gives us an orientation of $S^{n-1}_x$. It is easy 
to see 
that this orientation of $S_x$ does not depend of choice 
of the 
chart. Now, the pair (the orientation of $M$, the 
orientation of 
$S_x$) gives us an orientation of $STM$ which we fix 
forever. 

\begin{defin}\label{propaglaw} 
We define a {\it propagation law} on $M$ to be a smooth 
map 
$$ 
L: \overline{TM} \times \R \times \R \to \overline {TM} 
$$ 
(a time-dependent flow on $\overline{TM}$). Here $L 
({\mathbf u}, s, 
t)\in \overline {TM}$ should be thought of as the point 
that 
corresponds to the position and the velocity vector at 
moment 
$s+t$ of a perturbation whose position and the velocity 
vector at 
moment $s$ was ${\mathbf u}$. (We assume that a velocity 
of movement of a 
perturbation is either zero all the time or nonzero all 
the time.) 
Furthermore we assume that $L({\mathbf u},s,t)$ satisfies 
the following 
natural conditions: 
\begin{description} 
\item[a] $L({\mathbf u},s,0)={\mathbf u}$ for all 
${\mathbf u} \in 
\overline{TM}$; 
\item[b] $\forall s, 
t\in \R$ the map $L_{s,t}: \overline{TM} \to \overline{TM} 
$ defined as 
$L_{s,t}({\mathbf u})=L({\mathbf u},s, t)$ is a 
diffeomorphism; 
\item[c] $L({\mathbf u}, s, t_1+t_2)=L\bigl (L({\mathbf 
u}, s, t_1), s+t_1, 
t_2\bigr )$, 
$\forall s, t_1, t_2\in \R$; 
\item[d] $\forall {\mathbf u}\in \overline{TM}$ and 
$\forall s_0,t_0\in \R$ 
$$ 
\frac{d}{dt}p_T(L({\mathbf u}, s_0, t))\Big\vert_{t=t_0}=L 
({\mathbf u}, s_0, 
t_0). 
$$ 
\end{description} 
\end{defin} 

\begin{defin}\label{propag} 
A {\it propagation} is a quadruple $P=(L,x,T,V)$ where $L$ 
is a 
propagation law, $x\in M, T\in \R$ and $V:S^{n-1} 
_x\rightarrow 
\bigl ( T_xM\setminus s(x)\bigr )$ is a smooth section of 
the 
$\R^+$-bundle $\bigl(T_xM\setminus s(x)\bigr)\rightarrow 
S^{n-1}_x$. We fix an orientation preserving 
diffeomorphism 
$S^{n-1} \to S^{n-1}_x$ and further in the text regard $V$ 
as a 
mapping $V:S^{n-1}\rightarrow \bigl (T_xM\setminus s(x) 
\bigr )$. 
\end{defin} 

\m 
A propagation $P=(L,x, T,V)$ produces a wave front 
$W(t):S^{n-1}\rightarrow M, t \ge T,$ 
as follows. Informally speaking, we assume that at a 
moment of 
time $T$ something happens at a point $x\in M$ and the 
perturbation 
caused by this event starts to radiate 
from the point $x$ in all the directions according to a 
propagation law $L$ with the initial velocities of 
propagation in $T_xM$ 
described by $V$. Formally, for $t\geq T$ we define the 
front $W(t)$ to be the mapping 
$$ 
W(t):=p_T(L(V, T, t-T)):S^{n-1}\rightarrow M. 
$$ 
We put $\overline W(t)=L(V, T, t-T)$ and $\wt W(t) 
=\overline p 
\circ \overline W(t)$. In this case we also say that the 
wave 
front has {\it originated from the event $(x,T)$}. 
Initially a 
front of an event is a smooth embedded sphere (because of 
\ref{propaglaw}(d)), but generically it soon acquires 
double 
points, folds, cusps, swallow tails, and other complicated 
singularities. Generally, singular values of the front 
form a codimension two subset of $M$. 

\m 
We denote by $\eps _x: S^{n-1} \to STM$ any map of the 
form 
\begin{equation}\label{eps} 
\CD 
S^{n-1}@>h>> S^{n-1}_x\subset STM 
\endCD 
\end{equation} 
where $h$ is a map of degree 1. Clearly, the homotopy 
class of 
$\eps_x$ is well-defined and does not depend on $x$. 

\m 
Let $\mathcal S$ be the space of smooth maps $S^{n-1} 
\rightarrow STM$ 
that are homotopic to a map $\eps_x$ as in \eqref{eps}. 
Then $\mathcal S \times \mathcal S$ is the space of 
ordered 
pairs $(f_1, f_2)$ with $f_i \in \mathcal S$. 

Put $\Sigma$ to be the discriminant in $\mathcal S \times 
\mathcal S$, 
i.e. the subspace that consists of pairs $(f_1, f_2)$ such 
that 
there 
exist $y_1, y_2\in S^{n-1}$ with $f_1(y_1)=f_2 (y_2)$. 
(We do not include into $\Sigma$ the maps that are 
singular in the 
common sense but do not have double points between the two 
different spheres.) 

\begin{defin}\label{sigma0} 
We define $\Sigma_0$ to be a subset (stratum) of $\Sigma$ 
consisting 
of all 
the pairs $(f_1, f_2)$ such that there exists precisely 
one pair of 
points 
$y_1, y_2\in S^{n-1}$ such that: 
\begin{description} 
\item[a] $f_1(y_1)=f_2(y_2)$. And moreover this pair of 
points is such that: 

\item[b] $y_i$ is a regular point of $f_i, i=1,2$; 
\item[c] $(df_1)(T_{y_1})\cap (df_2)(T_{y_2})=0$. Here 
$df_i$ is the 
differential of $f_i$ and $T_{y_i}$ is the tangent space 
to $S^{n-1}$ 
at $y_i$. 
\end{description} 
\end{defin} 

\begin{con}\label{vector} 
Let $\rho:(a,b)\rightarrow \mathcal S\times \mathcal S$ be 
a path 
which intersects $\Sigma_0$ in a point $\rho(t_0)$. We 
also assume that 
$$ 
\rho(t_0-\delta, t_0+\delta)\cap \Sigma_0=\rho(t_0) 
$$ 
for $\delta$ small enough. We construct a vector ${\mathbf 
v}={\mathbf v}(\rho,t_0, \delta)$ as follows. 
We regard $\rho(t_0)$ as a pair $(f_1,f_2)\in \mathcal S 
\times \mathcal S$ 
and consider the points $y_1,y_2$ as in \ref{sigma0}. Set 
$z=f_1(y_1)=f_2(y_2)$. 
Choose a small $\delta>0$ and regard $\rho(t_0+\delta)$ as 
a pair 
$(g_1,g_2)\in \mathcal S \times \mathcal S$. Set 
$z_i=g_i(y_i), i=1,2$. 
Take a chart for $STM$ that contains $z$ and $z_i,\, i=1,2 
$ 
and set 
$$ 
{\mathbf v}(\rho,t_0,\delta) := 
\overrightarrow{zz_1}-\overrightarrow{zz_2}\in T_zSTM. 
$$ 
\end{con} 

\begin{defin}\label{transversal} 
Let $\rho:(a,b)\rightarrow \mathcal S\times \mathcal S$ be 
a path 
as in \ref{vector}. We say that $\rho$ {\it intersects 
$\Sigma_0$ 
transversally for $t=t_0$} if there exists $\delta_0>0$ 
such that 
\[ 
{\mathbf v}(\rho, t_0, \delta) \notin (df_1)(T_{y_1}S^{n- 
1})\oplus 
(df_2)(T_{y_2}S^{n-1}) \subset T_zSTM 
\] 
for all $\delta \in (0, \delta_0)$. 
\end{defin} 

It is easy to see that the concept of transversal 
intersection 
does not depend on the choice of the chart. 

\begin{defin}\label{genericpath} 
A path $\rho:(a,b)\rightarrow \mathcal S\times \mathcal 
S, -\infty 
\le a<b\le \infty$ is said to be {\it generic} if 
\begin{description} 
\item[a] $\rho(a,b)\cap \Sigma =\rho(a,b)\cap \Sigma_0$; 
\item[b] the set $J=\{t|\rho(t)\cap \Sigma_0\ne \emptyset\} 
\subset (a,b)$ 
is an isolated subset of $\R$; 
\item[c] the path $\rho$ intersects $\Sigma_0$ 
transversally for all $t\in J$. 
\end{description} 
\end{defin} 

As one can expect, every path can be turned into a generic 
one by 
a small deformation. We leave a proof to the reader. 

\m 
Let $P_1=(L_1,x_1,T_1, V_1)$ and $P_2=(L_2,x_2,T_2, V_2)$ 
be two propagations. They define mappings 
$r_i:\R\rightarrow 
\mathcal S$, $i=1, 2$ as follows. 
$$ 
r_i(t)= 
\begin{cases} 
\overline p \circ V_i & \text{ for }t\leq T_i,\\ 
\wt W_i(t) & \text{ for } t> T_i.\\ 
\end{cases} 
$$ 

\begin{defin}\label{genericpair} 
A pair of propagations $\{ P_1, P_2 \}$ is said to be {\it 
generic} 
if the path $r=(r_1, r_2):\R\rightarrow \mathcal S \times 
\mathcal 
S$ is generic and $r(T_i)\not \in \Sigma$, $i=1,2$. 
\end{defin} 

\begin{defin}\label{dang} 
Let $\{P_1, P_2 \}$ be a generic pair of propagations and 
let 
$r:\R \rightarrow \mathcal S \times \mathcal S$ be as 
above. Then 
a moment $t\in \R$ such that $r(t)\in \Sigma$ corresponds 
either 
to the baby-intersection or to the case where $t>\max(T_1, 
T_2)$ 
and there exists $y_1, y_2\in S^{n-1}$ such that 
$W_1(t)(y_1)=W_2(t)(y_2)=z$ and $\wt W_1(t)(y_1)=\wt 
W_2(t)(y_2)\in S^{n-1}_z$, i.e. to the case where there is 
a 
double point of the two fronts $W_1(t)$ and $W_2(t)$ at 
which the 
directions of the propagations of the two fronts coincide. 
Such a 
double point of two fronts is called {\em a point of 
dangerous 
intersection.} Notice that we do not exclude situations in 
which the 
two fronts are tangent, the so-called {\it Arnold's 
dangerous 
tangencies}, cf. \cite{Arnold}. 
\end{defin} 

\m 
For many pairs of propagations the dangerous intersection 
points 
do not occur. Such pairs of propagations are called {\it 
dangerous intersections 
free.} 
Now we describe a source of examples of such pairs. 

\begin{sourex}\label{source} 
Let $L:\overline{TM}\times \R \times \R\rightarrow 
\overline{TM}$ 
be a propagation law. Suppose that there exists a section 
$$ 
\wt s: STM \times \R \times \R \rightarrow \overline{TM} 
\times \R \times \R 
$$ 
of the map $\overline p \times 1 \times 1$ such that $\Im 
(\wt s)$ 
consists of the trajectories of $L$, i.e. if $({\mathbf 
u}, s_0, 0)\in \Im 
(\wt s)$ for some ${\mathbf u}\in \overline TM$ and $s_0 
\in \R$, 
then $L({\mathbf u}, s_0, t)\in \Im(\wt s)$, for every 
$t\in \R$. 

Let $P_1=(L, x_1, T_1, V_1)$ and $P_2=(L, x_2, T_2, V_2)$ 
be propagations such 
that $(\Im (V_1), T_1, 0)\subset \Im (\wt s\big |_{S_ 
{x_1}, T_1, 0})$, 
$(\Im (V_2), T_2, 0)\subset \Im (\wt s\big |_{S_{x_2}, 
T_2, 
0})$ and $r(T_i)\not \in \Sigma$, $i=1,2$. 

It is easy to see that {\it the pair $(P_1,P_2)$ is 
dangerous 
intersections free.} 
\end{sourex} 

\begin{ex}\label{metricpropagation}{\bf Propagations that 
are defined 
by a Riemannian metric.\/} An interesting class of 
examples comes 
from the propagation defined by the geodesics of a 
complete 
Riemannian metric $g$ on $M$. In this case $L({\mathbf u}, 
s, t)$ is just 
a point on $\overline {TM}$ that corresponds to a velocity 
vector 
at moment $s+t$ of a geodesic curve that had a velocity 
vector 
${\mathbf u}$ at the moment $s$. Thus the wave front $W(t) 
$ corresponding to a 
propagation $(L,x,T,V)$ can be described 
as $W(t):S^{m-1}\to M$ with $W(t)(y)=\exp_x\bigl ((t-T)V(y) 
\bigr), y\in 
S^{m-1}$, 
where $\exp_x:T_xM\to M$ is the exponent map corresponding 
to the 
Riemannian metric $g$. 

It is easy to see that in these examples if $V(S^{m-1})$ 
is a sphere of some 
radius $r$, then at 
every moment of time the velocity vectors of the points on 
the 
wave front are perpendicular to the image of the front. 
{\em Thus, in this case the dangerous intersections are 
precisely 
the dangerous tangencies.\/} 

Furthermore, if both $\Im V_1$ and $\Im V_2$ are spheres 
of the 
same radius $r$, {\em then the dangerous intersections (= 
dangerous 
tangencies) do not occur,} since spheres of radius $r$ in 
all the 
tangent spaces produce the section $\wt s$ described 
above. 
\end{ex}

\begin{ex}\label{nonisotropicpropagation}{\bf Propagation 
in a non-homogeneous 
and non-isotropic 
medium whose structure does not depend on time.\/} Assume 
that $M$ is a 
Riemannian manifold and $\mu: STM\rightarrow \overline {TM} 
$ is a 
smooth section of the corresponding $\R^+$-bundle such 
that $\Im 
(\mu\big |_{S_x})$ bounds a strictly convex domain in 
$T_xM$ for all $x\in M$. 
The radius vector from $s(x)$ to $\Im (\mu \big |_{S_x})$ 
in the given direction is the velocity vector of the 
distortion traveling 

in the direction. This information allows us 
to calculate for every smooth curve $\gamma: [t_1, t_2] 
\rightarrow 
M$ the total time $\tau(\gamma)$ needed for the distortion 
to travel 
along this curve. 

Assume $\Im (V_1)\subset \Im \mu$ and $\Im (V_2)\subset 
\Im \mu$ and that 
propagation occurs according to the Huygens principle, 
i.e. distortion 
travels along the extremal curves of the functional $\tau$ 
on the space 
of smooth curves on $M$. {\em It is clear that here we 
have a special case of 
the situation described in \ref{source}, and so the 
dangerous intersection 
points do not occur for such a pair of propagation.\/} On 
the other 
hand, if the propagation happens according to the Huygens 
principle then 
at every point $W(t)(x)=z$ the normal vector to the wave 
front is 
conjugate with respect to $\mu\big |_{S_z}$ to the 
direction of 
the extremal curve along which the information traveled to 
this 
point, 
Arnold~\cite{Arnoldbook}. In 
particular, in this case the dangerous tangencies do not 
occur 
under the wave fronts propagation, since they are the 
dangerous 
intersections. 
\end{ex}

\section{The causality relation invariant}\label{cri}

\m Recall that the standard sphere $S^{n-1}$ is assumed to 
be oriented. 
We say that a tangent frame $\mathfrak r$ to $S^{n-1}$ is 
positive if 
it gives us the standard orientation of $S^{n-1}$. 

\begin{defin}\label{signalk} 
Let $\rho$ be a path in $\mathcal S\times \mathcal S$ that 
intersects 
$\Sigma$ transversally in one point $\rho(t_0)\in \Sigma_0 
$. We 
associate a sign $\wt\sigma(\rho, t_0)$ to such a crossing 
as follows. 

We regard $\rho(t_0)$ as a pair 
$(f_1,f_2)\in \mathcal S \times \mathcal S$ and consider 
the points 
$y_1,y_2\in S^{n-1}$ such that $f_1(y_1)=f_2(y_2)$. 
Set $z=f_1(y_1)=f_2(y_2)$. 
Let $\mathfrak r_1$ and $\mathfrak r_2$ be frames which 
are 
tangent to $S^{n-1}$ at $y_1$ and $y_2$, respectively, 
and both are assumed to be positive. Consider the frame 
$$ 
\{ df_1(\mathfrak r_1), {\mathbf v}, df_2(\mathfrak r_2) 
\} 
$$ 
at $z\in STM$ where ${\mathbf v}$ is a vector described in 
\ref{vector}. We put $\wt\sigma(\rho, t_0)=1$ if this 
frame gives 
us the orientation of $STM$, otherwise we put $\wt\sigma 
(\rho, 
t_0)=-1$. 

Because of the transversality and condition (c) from \ref 
{sigma0}, 
the family $ \{ df_1(\mathfrak r_1), {\mathbf v}, df_2 
(\mathfrak 
r_2) \} $ is really a frame. 

Notice also that the vector ${\mathbf v}$ is not well- 
defined, but 
the above defined sign $\wt\sigma$ is. 
\end{defin} 

Clearly if we traverse the path $\rho$ in the opposite 
direction 
then the sign of the crossing changes. 

\begin{defin}\label{positive} 
Suppose that a front $W$ passes through a point $x\in M$ 
at the 
moment of time $t_0$ in such a way that the velocity 
vector 
${\mathbf v}_x$ of the front $W(t_0)$ at $x$ is transverse 
to 
$W(t_0)$, $W(t_0)$ restricted to a small neighborhood $U$ 
of 
$W^{-1}(t_0)(x)$ is an embedding, and $x$ has only one 
preimage 
under $W(t_0)$. 

Recall that the manifold $M$ is oriented. Let $o_x$ be the 
local 
orientation of $W(t_0)$ at $x$ (i.e. the orientation of 
the 
tangent plane $T_x$ to $W(t_0)$). We say that the local 
orientation $o_x$ is {\em positive}, and write $\sigma 
(W(t_0),x)=1$ if the pair ($o_x$, ${\mathbf v}_x$) gives 
us the orientation 
of $M$; otherwise we say that the local orientation of $W 
(t)$ at 
$x$ is {\em negative} and write $\sigma (W(t_0),x)=-1$. 
\end{defin} 

\begin{figure}[htbp] 
\begin{center} 
\epsfxsize 5cm 
\hepsffile{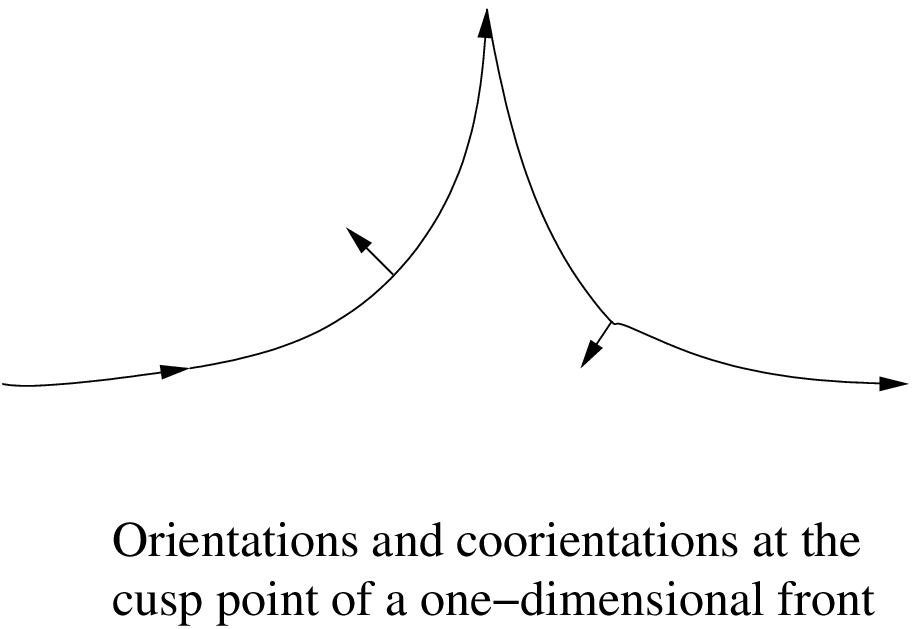} 
\end{center} 
\caption{}\label{orient.fig} 
\end{figure}

\m Notice that the same wave front $W(t)$ 
can contain two points $x$ and $y$ such $o_x$ is positive 
orientation 
while $o_y$ is the negative one, 
see Figure~\ref{orient.fig}. 

\m Consider a generic pair $(P_1,P_2)$ of propagations 
$P_1=(L_1,x_1, V_1,T_1)$ and $P_2=(L_2,x_2, V_2,T_2)$. 
{\bf In the rest of 
Section~\ref{cri} we 
assume that 
$T_1\leq T_2$. The case where $T_1>T_2$ is treated in a 
similar way.\/} 

Let $t> T_2$ be a generic moment of time, i.e. the one 
at which dangerous intersections do not occur. 

Let $c_i, i\in I \subset \N$ where 
$T_2<c_i< t$ be moments of time when dangerous 
intersections did occur. 

\begin{defin}\label{signcri} 
We define $\sigma(W_1(c_i), W_2(c_i))$ as the sign $\wt 
\sigma$ 
of the corresponding passage of $\Sigma_0$. 
\end{defin} 

Notice that 
$\sigma(W_1(c_i), W_2(c_i))$ is symmetric if $n$ is even 
and skew-symmetric if $n$ is odd. 

\m 
Let $p_j,j\in J\subset \N$ be the moments of time when the 
front 
$W_1$ passed through the point $x_2$ before the front $W_2 
$ 
originated. (Notice that $p_j<T_2$ and that 
$\sigma(W_1(p_j),x_2))$ is well-defined since the pair of 
propagations is generic.) A straightforward verification 
shows 
that 
$$ 
\sigma(W_1(p_j),x_2)=\wt\sigma(\rho, t_0). 
$$ 
where $\rho(t)=(\wt W_1(t), \eps_{x_2}), t\in (p_j-\delta, 
p_j+\delta)$.

\begin{defin} 
We set 
\begin{equation*} 
\cri (W_1(t), W_2(t))=\sum_{i\in 
I}\sigma(W_1(c_i), W_2(c_i))+\sum_{j\in J}\sigma(W_1 
(p_j),x_2)\in \Z 
\end{equation*} 
and call it the {\it causality relation invariant} for the 
fronts $W_1(t)$ and 
$W_2(t)$ at a given moment of time $t$. (If in fact 
$T_1>T_2$, then the second 
sum should be $\sum_{k\in K}(-1)^{\dim M}\sigma (W_2 
(q_k),x_1)$, where $q_k,k\in 
K\subset \N$ are the moments of time when the front 
$W_2$ passed through the point $x_1$ before the front $W_1 
$ 
originated. One can easily verify that $(-1)^{\dim M} 
\sigma (W_2(q_k),x_1)$ 
coincides with the sign of the corresponding crossing of 
$\Sigma_0$ by the path 
$r=(r_1,r_2)$.) 
\end{defin}

So, if $\cri (W_1(t), W_2(t))=k\neq 0$ then the sum of the 
number of 
baby-intersections and of the number of dangerous 
intersections is 
at least $|k|$. 
This probably could be interpreted as the quantity 
that measures either how much faster the first front is 
than the 
second so that they could become dangerously intersected; 
or how 
many times the first front did pass through the source of 
the 
second front before the second front originated. 

\m Now we consider an important special case. 
{\it If the above pair $(P_1,P_2)$ of propagations is 
dangerous intersections 
free, then\/} 
\[ 
\cri (W_1(t),W_2(t))=\sum_{j\in J} \sigma (W_1(p_j),x_2). 
\] 

We saw examples of such propagations in 
\ref{source},~\ref{metricpropagation},~\ref 
{nonisotropicpropagation}. 
It is easy to see that in this case $\cri (W_1(t), W_2(t)) 
$ 
does not depend 
on $t$ provided $t>T_2$, and thus it is invariant 
under the propagation. In particular, if $\cri (W_1(t), W_2 
(t))$ 
is non-zero, then we know for a fact that the 
perturbation caused by the first signal has reached the 
source point of the second signal before the second signal 
originated. Moreover, if $\cri (W_1 (t), W_2(t))=k\neq 0$, 
then we can say for sure that the first wave 
front has passed through the source point of the second 
front 
at least $|k|$ times before the second signal originated. 
(Of course it could be that it did pass more 
times, because it could have passed $k+l$ times with a 
positive sign and $l$ times with a negative sign.) Thus, 
$W_1$ and $W_2$ 
can be thought of as being causally related.

In case $\cri(W_1, W_2)=0$ we can not make any conclusions 
on 
causality relation between $W_1$ and $W_2$, since it could 
mean that 

\begin{description} 
\item[1] neither of the fronts passed through the origin 
of the other before the 
other was born, and thus the fronts are not causally 
related; or 
\item[2] the earlier-born front did pass through the 
origin of the other front 
$2k\neq 0$ times before the other front was born with $k$ 
of these 
passages being positive and $k$ of the passages being 
negative. In this case 
the fronts are causally related. 
\end{description}

However, for some Riemannian metrics $g$ the second 
situation could not occur, 
as the following example shows. 

\begin{ex}[when $\cri(W_1(t), W(t)=0$ if and only if the 
two fronts are causally unrelated]\label{Friedmann} 
Consider the case of the propagation law given by a 
complete Riemannian metric, 
see~\ref{metricpropagation}, and assume that the 
propagations 
$P_1=(L,x_1, T_1, V_1)$ and 
$P_2=(L,x_2, T_2, V_2)$ be such that $V_1(S^{m-1})$ and 
$V_2(S^{m-1})$ are 
spheres of the same radius $r$ in $T_{x_1}M$ and $T_{x_2} 
M$, respectively. As it 
was discussed in~\ref{metricpropagation}, in this case 
dangerous intersection 
points between the two fronts are the dangerous tangency 
points and they do not 
occur during the propagation of $W_1$ and of $W_2$. 

Suppose now that $M$ is a manifold of non-positive 
sectional curvature. Then, by 
the Hadamard Theorem, see for example~\cite{docarmo}, the 
exponential map 
$\exp_x: T_xM \to M$ is the universal covering map for 
every $x\in M$. Thus, 
the front $W_i(t): \exp_{x_i}((t-T_i)V): S^{m-1} \to M, 
i=1,2,$ is always 
an immersion, and the local orientation of a wave front 
passing through 
every point $x$ is positive. In particular, all the 
passages of $W_1$ through 
$x_2$ 
are positive. Thus these passages can not cancel each 
other in the definition of 
$\cri(W_1, W_2)$. And we get that for such cases $\cri 
(W_1, W_2)\neq 0$ 
if and only if the fronts are causally related. 

Propagation laws given by a Riemannian metric of constant 
sectional curvature 
with 
$V_1$ and $V_2$ being spheres of unit radius are extremely 
important in 
cosmology 
and they are known as {\it Friedmann Cosmology models}, 
see for 
example~\cite{Frankel},~\cite{Peebles},~\cite 
{CornishWeeks},~\cite{Sormani}. 
Thus we conclude that for Friedmann cosmology models given 
by metrics of 
constant non-positive sectional curvature $\cri(W_1, W_2) 
\neq 0$ if and only if 
$W_1$ and $W_2$ are causally unrelated. 
\end{ex}

\section{Homotopy properties of maps to $STM$}\label{hom}

\begin{defin}\label{spec} 
Given a map $\alpha: S^1\times S^{n-1}\rightarrow STM$, we 
say that $\alpha$ is {\it special} if 
$\alpha\big |_{* \times S^{n-1}}$ has the form $\eps_x$ 
for 
some $x\in M$, see~\eqref{eps}. Here $*\in S^1$ is the 
base point. 
\end{defin} 

\begin{defin} 
Given an $n$-dimensional manifold $N$ and a map $\beta: 
N\rightarrow STM$, we 
define $d(\beta)$ to be the degree of the map 
\begin{equation*} 
\pr\circ\beta: N\rightarrow M 
\end{equation*} 
\end{defin}

\begin{lem}\label{degreesphere} 
Let $\alpha:S^1\times S^{n-1}\rightarrow STM$ be a special 
map. 
Then there exists a map $\beta : S^n\rightarrow STM$ such 
that 
$d(\beta) = d(\alpha)$. 
\end{lem} 

\pp We regard $S^{n-1}$ as a pointed space. Consider a map 
$\wt \alpha: 
S^1\times S^{n-1} \rightarrow STM$ such 
that: 
\begin{description} 
\item[1]$\wt \alpha\big |_{* \times S^{n-1}}=\alpha\big |_ 
{* \times S^{n-1}}$, 
\item[2] $\wt \alpha \big |_{S^1\times *}=\alpha \big |_ 
{S^1 \times *}$, 

\item[3] $\wt \alpha |_{t \times S^{n-1}}=\eps\big |_{\wt 
\alpha (t\times *)}$. 
\end{description} 
We regard $S^1\times S^{n-1}$ as the $CW$-complex with 
four 
cells $e^0, 
e^1,e^{n-1},e^n$, $\dim e^k=k$. 
It is easy to see that the maps $\wt \alpha$ and $\alpha$ 
coincide 
on the $(n-1)$-skeleton. Thus, the maps $\alpha$ and $\wt 
\alpha$ (restricted to 
the $n$-cell) together yield a map $\beta: 
S^n\rightarrow STM$. Clearly $d(\wt \alpha)=0$, and 
therefore 
$d(\beta)=d(\alpha)$. 
\qed 

\begin{lem}\label{Eulerclass} 
Suppose that there exists a map $\beta: S^n\rightarrow 
STM$ 
with $d 
(\beta)\neq 0$. Then the 
Euler class $\chi \in H^n(M)$ of the tangent bundle 
$TM\rightarrow M$ 
is zero. 
\end{lem}

\pp We set $f=\pr\circ\beta: S^n \to M$. Clearly, $M$ is 
closed and 
$H^n(M)=\Z$ because $d(\beta)\ne 0$. 
So, again since $d(\beta)\ne 0$, we conclude $f^*\chi \ne 
0$ whenever $\chi\ne 0$. Now the result follows because 
$f^*\chi$ is the obstruction to the lifting of $f$ to 
$STM$, 
while $\beta$ is such a lifting of $f$. 
\qed 

\begin{lem}\label{rationalhomologysphere} 
Let $M^n$ be an oriented manifold and $\beta: 
S^{n}\rightarrow 
STM$ be a map with 
$d(\beta)\neq 0$, then $M$ is a closed manifold which is a 
rational homology 
sphere. 
\end{lem} 

\pp 
We set~$f=\pr\circ\beta:S^{n}\to M$ and $d=d(\beta)$. 
Clearly $M$ is 
closed because 
$d(\beta)\ne 0$. Let 
$f_!:H_*(M) \to H_*(S^{n})$ be the transfer map, see 
e.g. \cite[V.2.12]{Rudyak}. 
Since $f_*(f^*y\cap x)=y\cap f_*x$ for all $x\in H_*(S^{n}) 
$ 
and $y\in H^*(M)$, we conclude that $f_*f_!(z)=dz$, for 
all $z\in H_*(M)$. In 
particular, since $H_i(S^n)=0,$ for $0<i<n$, then $dH_i(M) 
=0,$ for $0<i<n$. 
Thus $H_i(M;\Q)=0$ for $0<i<n$ and $M$ is a rational 
homology sphere. 
\qed 

\begin{cor}\label{evencase} 
If $M^{2k}$ is an even-dimensional oriented manifold, then 
$d(\beta)=0$ for every $\beta: S^n\rightarrow STM$. 
\end{cor} 

\pp By Lemma~\ref{rationalhomologysphere} we get that if $d 
(\beta)\neq 0$ for 
some $\beta:S^{2k}\rightarrow STM$, then $H_i(M;\Q)=0$ 
for $0<i<2k$ and $H_0(M,\Q)=H_{2k}(M,\Q)=\Q$. Thus the 
Euler characteristic of $M$ is 2 and the Euler class 
of the tangent bundle $TM\to M$ is non-zero (in fact, $\pm 
2$). 
This contradicts to the statement of Lemma~\ref 
{Eulerclass}. 
\qed 

\m 
Let $\deg: \pi_n(M^n) \to \Z$ be the degree homomorphism, 
i.e., the homomorphism 
which assigns the degree $\deg f$ to the homotopy class of 
a map $f: S^{n} \to M$. 
(In fact, it coincides with the Hurewicz homomorphism $h: 
\pi_n(M) \to H_n(M)$ 
for $M$ closed and is zero for $M$ non-closed.) 

\begin{defin}\label{a} 
Given a connected oriented manifold $M^n$, we define an 
Abelian group $\A(M)$ 
and a homomorphism $q=q_M:\Z\to \A(M)$ as follows. If $n$ 
is even, then 
$\A(M)=\Z$ and $q=1_{\Z}$. If $n$ is odd, then $\A(M)$ is 
the cokernel of the 
degree homomorphism $\deg: \pi_n(M) \to \Z$ and $q: \Z \to 
\A(M)$ is the 
canonical epimorphism. 
\end{defin}

\begin{prop}\label{usefulprop} 
Let $M^n$ be an odd-dimensional manifold as in 
Definition $\ref{a}$. 
Then the following holds: 
\par{\rm (i)} If the universal covering space of $M$ is a 
non-compact manifold, 
then $\A(M)=\Z$. 
\par{\rm (ii)} If $M$ admits a complete Riemannian metric 
of non-positive 

sectional curvature, then $\A(M)=\Z$. 
\par{\rm (iii)} If $M$ is not a rational homology sphere, 
then $\A(M)=\Z$. 
\par{\rm (iv)} If $\pi_1(M)$ is infinite, then $\A(M)=\Z$. 
\par{\rm (v)} If $\pi_1(M)$ is a finite group of order 
$k$, then $A(M)=\Z/m\Z$ 
where $k|m$ {\rm (}the case $m=0$. i.e. $A=\Z$ is also 
possible{\rm )}. 
\par{\rm (vi)} If $M$ is a closed manifold with $\A(M)=0$, 
then $M$ is a 
homotopy sphere. 

\end{prop}

\pp (i) This follows because every map $S^n \to M$ passes 
through the universal covering $\wt M \to M$ while $H_n 
(\wt M)=0$. 
Therefore the Hurewicz homomorphism $h:\pi_n (M) \to H_n(M) 
$ is trivial. 

\par(ii) The Hadamard Theorem says that the universal 
covering of such $M^n$ is 
diffeomorphic to $\R^n$, see for example~\cite{docarmo}, 
and the statement 
follows from (i). 

\par(iii) This follows immediately from Lemma~\ref 
{rationalhomologysphere}. 

\par (iv) This follows from (i) since the universal 
covering space of $M$ is 
non-compact. 

\par (v) This follows because every map $S^n \to M$ passes 
through the universal covering map $p: \wt M \to M$ which 
is of degree $k$.

\par(vi) If $\A(M)=0$ then there exists a map $S^n \to M$ 
of 
degree 1. Since every map of degree 1 induces epimorphism 
of 
fundamental groups and homology groups, we conclude that 
$M$ is a 
homotopy sphere. \qed

\section{The affine linking invariant $\lk$ as a reduction 
of $\cri$}\label{lk}

\begin{defin}\label{sigma1} 
We define $\Sigma_1$ to be the subset (stratum) of 
$\Sigma$ consisting 
of all the pairs $(f_1, f_2)$ such that there exists 
precisely two 
pairs of points $y_1, y_2\in S^{n-1}$ as in \ref{sigma0}. 
Here we assume 
that the two double points of the image are distinct. 
\end{defin} 

Notice that $\Sigma_i$ is a stratum of codimension $i$ 
in $\Sigma$. In particular, a generic path in $\mathcal 
S\times 
\mathcal S$ intersects $\Sigma_0$ in a finite number of 
points, 
and a generic disk in $\mathcal S\times \mathcal S$ 
intersects 
$\Sigma_1$ in a finite number of points and does not 
contain singular 
mappings that are not in $\Sigma_0\cup \Sigma_1$. 

\m A generic path 
$\gamma :[0,1]\rightarrow \mathcal S \times \mathcal S$ 
that connects two points in $\mathcal S\times \mathcal 
S\setminus \Sigma$ intersects $\Sigma_0$ in finitely many 
points $\gamma(t_j), 
j\in J\subset N$ and all the intersection points are of 
the 
types described 
in~\ref{signalk}. Put 
\begin{equation}\Delta_{\lk}(\gamma)=\sum_{j\in J} \sigma 
(\gamma,t_j)\in \Z. 
\end{equation} 

\m We let $A=\{(x,y)\in \R^2\bigm| x^2+y^2\le 1\}$, 
$B_1=\{(x,y)\in A \bigm| 
xy=0\}$, $B_2=\{(x,y)\in A \bigm| 
x=0\}$, $B_3=\{(0,0)\}$, $B_4=\emptyset$. 

We define a {\it regular disk} in $\mathcal S\times 
\mathcal S$ 
as an embedded disk $D$ such that the triple $(D, D \cap 
\Sigma_0, D \cap \Sigma_1)$ 
is homeomorphic to a triple $(A, B, C), A\supset B \supset 
C$ where $B$ is one 
of $B_i$'s and $C$ is a (possibly empty) subset of $B_3$. 

\begin{lem}\label{small} 
Let $\beta$ be a generic loop that bounds a regular disk 
in 
$\mathcal S\times \mathcal 
S$. Then $\Delta_{\lk}(\beta)=0$. 
\end{lem} 

\pp Straightforward. 
\qed 

\begin{lem}\label{disk} 
Let $\beta$ be a generic loop that bounds a disk in 
$\mathcal S\times \mathcal 
S$. Then $\Delta_{\lk}(\beta)=0$. 
\end{lem} 

\pp Without loss of generality we can (using a small 
deformation 
of the disk) assume that the disk is the union of regular 
ones, cf. 
Arnold~\cite {Arnold},~\cite{Arnoldcurves}. Now the proof 
follows 
from Lemma \ref{small}. 
\qed 

\m Notice that $\mathcal S\times \mathcal S$ is path 
connected. 

\begin{cor}\label{homomorphism} 
The invariant $\Delta_{\lk}$ induces a well-defined 
homomorphism 
$\Delta_{\lk}:\pi_1 (\mathcal S\times \mathcal 
S, *)\rightarrow \Z$. 
\end{cor} 

\pp Since every element of $\pi_1(\mathcal S\times 
\mathcal 
S, *)$ can be represented by a generic loop, the proof 
follows 
from Lemma \ref{disk}. 
\qed 

\m 
Let $x_1, x_2$ be two distinct points of $M$. Let $\alpha: 
S^1\times 
S^{n-1}\rightarrow STM$ be a special map (see Definition 
\ref{spec}) such that 
the composition 
\begin{equation*} 
S^{n-1} \subset S^1\times S^{n-1}\xrightarrow{\alpha} STM 
\end{equation*} 
has the form $\eps_{x_1}$, and let $e:S^1\times S^{n-1} 
\rightarrow STM$ be the map of the form 
\begin{equation*} 
S^1\times S^{n-1}\xrightarrow{\text{proj}} S^{n-1} 
\xrightarrow{\eps_{x_2}} STM 
\end{equation*} 
Then $(\alpha,e)$ is a loop in $(\mathcal S \times 
\mathcal 
S, *)$. 

\begin{lem}\label{degreeDelta} 
$\Delta_{\lk}[(\alpha,e)]=d(\alpha)$. 
\end{lem} 

\pp Notice that $\Delta_{\lk}[(\alpha,e)]$ is the 
intersection index of the 
cycles $\alpha(S^1\times S^{n-1})$ and $S^{n-1}_{x_2}$. 
This index coincides 
with the degree of the map $\pr\circ \alpha$ because the 
last one is equal to 
the algebraic number of the preimages of $x_2$. 
\qed 

\begin{defin} 
Take two different points $x_1, x_2\in M$ and consider the 
point 
$*=(\eps_{x_1},\eps_{x_2})\in\mathcal S\times \mathcal 
S\setminus 
\Sigma$, take an arbitrary point $f= 
(f^1_1,f^1_2)\in 
\mathcal S\times \mathcal S \setminus \Sigma$ and choose a 
generic path $\gamma$ 
going from $*$ to $f$. We set 
\begin{equation*} 
\lk(f)=q(\Delta_{\lk}(\gamma))\in \A(M) 
\end{equation*} 
and call $\lk$ the {\it affine linking invariant}. Here 
$q$ 
is the epimorphism from Definition~\ref{a}. 
\end{defin}

\begin{thm}\label{Main} 
The function $\lk:\pi_0(\mathcal S\times \mathcal S 
\setminus \Sigma) 
\rightarrow \A(M)$ is well-defined and increases by $1\in\A 
(M)$ under the 
positive transverse passage through the stratum $\Sigma_0 
$. 
\end{thm} 

\begin{rem} 
It is easy to see that $\lk$ does not depend on the choice 
of the pair 
$(x_1, x_2)$ used to define it. This follows, since all 
the pairs 
$(\eps_{x_1},\eps_{x_2}), x_1, x_2\in M$ are isotopic. 

However, if one chooses a basepoint $\tilde *$ in the 
construction of $\lk$ that is not two fibers over two 
distinct points, then the $\lk$-type invariant obtained 
this way would be different from $\lk$ by an additive 
constant $\lk(\tilde *)-\lk(*)$. This ambiguity is similar 
to the ambiguity in the choice of the zero vector in an 
affine vector space, and it is the reason for the 
adjective ``affine'' in the name of the invariants. The 
general theory of affine linking invariants is discussed 
in our work~\cite{ChernovRudyakGeneralLinking}. 
\end{rem}

{\em Proof\/} of Theorem~\ref{Main}.
To show that $\lk$ is well-defined we must verify that the 
definition 
is independent on the choice of the path $\gamma$ that 
goes 
from 
$*$ to $f$. This is the same as to show that $q(\Delta_ 
{\lk} 
(\varphi))=0$ for 
every closed generic loop $\varphi$ at $*$. 

Since $\Delta_{\lk}:\pi_1(\mathcal S\times \mathcal S, *) 
\rightarrow \Z$ is a 
homomorphism, it suffices to prove that $q(\Delta_{\lk} 
(\varphi))=0$ for all 
generators $\varphi$ of $\pi_1 (\mathcal S\times \mathcal 
S, *)$. 

The classes $[(\alpha,e)]$ and 
similar classes 
$[(e,\alpha)]$ generate the group $\pi_1(\mathcal S\times 
\mathcal S, *)$. By 
Lemma \ref{degreeDelta} we have 
\[ 
\Delta_{\lk}[(\alpha,e)]=d(\alpha). 
\] 
Clearly, $d(\alpha)=0$ if $M$ is an non-closed manifold. 
So, we 
assume $M$ to be 
closed. Now, for $n$ even $d(\alpha)=0$ by Corollary \ref 
{evencase}, while for 
$n$ odd $q(d(\alpha))=0$ by Lemma \ref{degreesphere}. 
\qed

\m 
Let $(P_1, P_2)$ be a 
generic pair of propagations, and let $t$ be a moment of 
time 
when dangerous intersection do not occur. 
Let $q: \Z \to \A(M)$ be the epimorphism described in 
Definition~\ref{a}. 

\begin{thm}\label{main} 
The invariants $\cri$ and $\lk$ are related as follows: 
$$ 
q\bigl(\cri(W_1(t), W_2(t)\bigr)=\lk\bigl(\wt W_1(t),\wt 
W_2(t)\bigr). 
$$ 
\end{thm} 

\pp The Theorem follows because the signs defined for 
dangerous 
intersections (as well as for baby-intersections) are 
exactly the sign 
of the corresponding crossings of $\Sigma_0$. 
\qed 

\m 
The invariant $\lk$ works especially nice for even 
dimensional manifolds, since in these cases $\A(M)=\Z$. 

\m If the pair of propagations is dangerous intersections 
free 
then the invariant $\lk(\wt W_1(t), \wt W_2(t))\in\A$ 
gives us the 
number of times the earlier-born front had passed through 
the birth point of the 
other front before the other front was born. 
It is also easy to see that if the other wave front did 
not appear yet, then 
$\lk (\wt W_1(t), \wt W_2(t))$ is the number of times the 
earlier-born 
front passed through the (future) birth point of the other 
front by the time 
moment $t$. 

\m We illustrate our general results with the following 
example. 

\begin{thm}\label{Friedmannlinking} Assume that the law of 
propagation $L$ is 
given by a complete Riemannian metric as in $\ref 
{metricpropagation}$, of 
non-positive sectional curvature on 
$M$ and the wave fronts $W_i(t), 
i=1,2,$ correspond respectively to 
propagations $P_i=(L, x_i, T_i, V_i), i=1,2,$ with $V_i:S^ 
{n-1}\to T_{x_i}M, 
i=1,2,$ being spheres of the same radius. Then $W_1$ and 
$W_2$ are causally 
related if and only $\lk (\wt W_1(t), \wt W_2(t))\neq 0$. 
\end{thm} 

\pp As it was explained in~\ref{Friedmann} for such 
propagations $\cri(W_1, 
W_2)\neq 0$ if and only if $W_1$ and $W_2$ are causally 
related. Notice that 
$\A(M)=\Z$ if $M$ is even-dimensional by definition of $\A 
(M)$, while $\A(M)=\Z$ 
for 
$M$ odd-dimensional by \propref{usefulprop}. Now the 
result follows from 
\theoref{main}. 
\qed 

\m This Theorem, in particular, says that $\lk(W_1(t), W_2 
(t))$ allows one to 
detect always whether $W_1$ and $W_2$ are causally related 
or not in Friedmann 
models based 
on metrics of constant non-positive sectional curvature.

\section{Causality relation invariant in the case of the 
propagation 
according to Riemannian metrics.}\label{mod} 

As it was noticed in~\ref{metricpropagation}, if a 
propagation 
happens according to a complete Riemannian metric and $\Im 
V$ is a 
sphere, then the velocities of the points of the front are 
always 
orthogonal to the front. So, if each of the two 
propagations happens 
according to a complete Riemannian metric, then dangerous 
tangency 
points and the dangerous intersection points are the same 
thing. 
In this section we deal only with this case. Namely, we 
provide an 
explicit way of calculation of the invariant $\cri$. We 
need some preliminaries. 

\m Let $W$ be a wave front, and let $x\in \Im W(t)$ be a 
non-singular point of $W(t)$. For sake of simplicity we 
denote 
$T_x(\Im W(t))$ just by $T$. Let $O$ be a small 
neighborhood of $x$ 
in $M$, and let $U=O\cap \Im W(t)$. Without loss of 
generality we 
can and shall assume that the injectivity radius is big 
enough 
($\ge 3$) for all points of $O$. 

The Riemannian metric $g$ on $M$ produces a unique 
symmetric 
connection on $M$. So, for every $a\in O$, the parallel 
transport 
along the geodesic segment (connecting $x$ and $a$) gives 
us an 
isomorphism 
\begin{equation}\label{parallel} 
{\mathbf \tau}_a: T_aM \to T_xM 
\end{equation} 
Furthermore, we can regard every sphere $S_a\in STM, a\in 
M$ as 
the unit sphere in $T_aM$, and so $STM$ can be regarded as 
the 
total space of the unit sphere subbundle of $TM$. Since 
the 
connection respects the Riemannian metric, we conclude 
that 
${\bf \tau}_a(S_a)=S_x$. 

\begin{defin}\label{setup} 
(a) We define 
$$ 
\pi: \pr^{-1}(O) \to S_x 
$$ 
as follows. A point $z\in \pr^{-1}(O)$ is a pair $(a,\xi)$ 
with 
$a=\pr(z)$ and $\xi \in S_a$, and we set $\pi(z)={\bf\tau} 
_a(\xi)$ 
with ${\mathbf\tau}_a$ as in \eqref{parallel}. 

(b) Let ${\mathbf n}_u$ be the unit normal vector to $U$ 
at $u$ that points to 
the direction of propagation of the front. We define the 
Gauss map $G=G_W: U \to 
S_x$ by setting 
$G(u)={\bf\tau}_u({\mathbf n}_u)$. 

(c) Given $u\in U$, let $\ell(u)={\mathbf n}_u\in \Im \wt 
W$. In this way we get 
a map $\ell: U \to \Im\wt W \subset STM$. We 
set $z=\ell(x)$. Given $\e \in T$, we set 
$$ 
\e^W:=d\ell(\e),\quad \e^W\in T_z\Im \wt W \subset T_zSTM. 
$$ 
It is clear that $d\ell: T \to T_z\Im\wt W$ is an 
isomorphism. 

(d) Let $z\in STM$ be the point described in (c). We 
define the 
{horizontal section} $H: O \to STM$ of $\pr$ by setting 
$$ 
H(a)={\bf\tau}^{-1}_a(z)\in S_a\subset STM. 
$$ 
Furthermore, given ${\mathbf w} \in T_aM, a\in O$, we set 
$$ 
{\mathbf w}^H=dH({\mathbf w})\in T_{H(a)}STM. 
$$ 
Clearly, ${\mathbf w}^H$ can be characterized by the 
properties 
$$ 
(d\pr)({\mathbf w}^H)= {\mathbf w}, \quad d\pi({\mathbf 
w}^H)={\mathbf 0}. 
$$

(e) Let $z\in STM$ be the point described in (c). Given 
${\mathbf 
w} \in T_xM$, we define ${\mathbf w}^S\in T_zT_xM$ as 
follows. We 
regard $z$ as the vector ${\mathbf z}\in T_xM$. 
Furthermore, we 
regard $T_xM$ as the affine space $T^{\aff}$ over the 
vector space 
$T_xM$ and consider the parallel shift 
$$ 
P_{\mathbf z}:T^{\aff} \to T^{\aff}, \quad a \mapsto a + 
{\mathbf 
z}. 
$$ 
Let $o\in T^{\aff}$ correspond to the origin of the vector 
space 
$T_xM$. Using the obvious identification $T_xM=T_oT^{\aff} 
$, we 
regard ${\mathbf w}$ as the tangent vector ${\mathbf w} 
_o\in 
T_oT^{\aff}$, and we set 
$$ 
{\mathbf w}^S=dP_{\mathbf z}({\mathbf w}_o)\in 
T_zT^{\aff}=T_zT_xM. 
$$ 
Notice that if $\e\in T$ then $\e^S \in T_zS_x$. (This is 
where 
the notation comes from: $\e^S$ is the spherical lifting 
of $\e$.) 
\end{defin} 

\begin{lem}\label{eh} 
For every $\e \in T$ we have $\e^W-\e^H=dG(\e)$. 
\end{lem} 

\pp First, notice that $G=\pi\circ \ell: U \to S_x$. 
So, $d\pi(\e^W)=dG(\e)$. Now, 
$$ 
(d\pr)(\e^W - dG(\e))=\e -{\mathbf 0}=\e 
$$ 
and 
$$ 
d\pi(\e^W - dG(\e))=dG(\e)-dG(\e)={\mathbf 0}. 
$$ 
Thus, $\e^W - dG(\e)=\e^H$. \qed 

\begin{prop}\label{nabla} 
Let $\e \in T$, and let ${\mathbf n}$ be the normal vector 
field 
to $U$ in $M$. Then $dG(\e)=(\nabla_{\e}{\mathbf n})^S.$ 
\end{prop} 

Here $\nabla$ denotes the covariant differentiation 
operation on $M$. 

\pp Let $\gamma: (-\delta, \delta) \to U$ be a curve with 
$\dot 
\gamma(0)=\e$. We define the curve $\zeta: (-\delta, 
\delta) \to 
S_x$ by setting $\zeta(t)$ to be the end of the vector 
${\bf\tau}_{\gamma(t)}{\mathbf n}_{\gamma(t)}$. Since 
$$ 
\nabla_{\e}{\mathbf n}=\frac 
d{dt}\left({\bf\tau}_{\gamma(t)}{\mathbf 
n}_{\gamma(t)}-{\mathbf n}_x\right)\Big\vert_{t=0}\, 
$$ 
we conclude that $\dot\zeta(0)= (\nabla_{\e}{\mathbf n}) 
^S$. On 
the other hand, $G\circ \gamma=\zeta$, and thus 
$$ 
dG(\e)=dG(\dot\gamma(0))= 
\dot{(G\circ\gamma)}(0)=\dot\zeta(0)=(\nabla_{\e}{\mathbf 
n})^S. 
$$ 
\qed 

\begin{cor}\label{covar} 
$\e^W-\e^H= (\nabla_{\e}{\mathbf n})^S$. 
\end{cor} 

\pp This is the direct consequence of \ref{eh} and \ref 
{nabla}. 
\qed 

\m 
Consider the Weingarten operator 
\begin{equation}\label{weingarten} 
A=A_W: T \to T, \quad A(\e)=\nabla_{\e}{\mathbf n}. 
\end{equation} 
The Corollary \ref{covar} can now be written as follows: 
\begin{equation}\label{oper} 
\e^W-\e^H= (A\e)^S. 
\end{equation} 

\m 
Now let $W_1$ and $W_2$ be two wave fronts, and let $x\in 
M$ be 
a point of dangerous tangency of $W_1(t)$ and $W_2(t)$. We 
assume 
that the corresponding pair of propagations is generic. 
Again, we 
denote by $T$ the common tangent plane $T_xW_i(t), i=1,2$. 
Let 
$A_i:=A_{W_i}: T \to T, i=1,2$ be the Weingarten operators 
considered in \eqref{weingarten}. We set $B=A_1 - A_2$. It 
is well 
known that each $A_i$ is a self-adjoint operator, \cite 
[Ch. 
7]{Spivak}, and therefore $B$ is. Let $k_1, \ldots, k_{n-1} 
$ be 
the eigenvalues (with multiplicities) of $B$. 

\begin{prop}\label{kerb} 
$\Ker B=0$, and so $k_i\ne 0$ for all $i$. 
\end{prop} 

\pp Let $\e\ne {\mathbf 0}$ be a vector with $B\e={\mathbf 
0}$. Then 
$$ 
\e^{W_1}-\e^{W_2}=(\e^{W_1}-\e^H)-(\e^{W_2}-\e^H)= 
(A_1\e-A_2\e)^S=(B\e)^S={\mathbf 0}. 
$$ 
i.e. $T_x\Im \wt W_1\cap T_x\Im \wt W_2\ne 0$. But this is 
impossible because 
the pair of propagations is assumed to be generic 
(see conditions \ref{genericpath}(b) and \ref{sigma0}(c)). 
\qed 

\m In particular, $\det B=k_1\cdots k_{n-1}\ne0$. 

\begin{defin}[{\bf Alternative definition of $\sigma(W_1 
(t), W_2(t))$}] 

We put $\eps(W_1(t), W_2(t)=1$ if both fronts have the 
same 
local orientations at $x$ (as defined in \ref{positive}) 
and 
$\eps(W_1(t), W_2(t))=-1$ if the fronts have opposite 
local 
orientations. Now we set 
\begin{equation*} 
\wh \sigma(W_1(t), W_2(t))=\eps(W_1(t), W_2(t))\sign(\det 
B) 
\sign(|{\mathbf v}_1|-|{\mathbf v}_2|) 
\end{equation*} 
where ${\mathbf v}_i$ is the velocity vector of $W_i(t)$ 
at $x$. 
\end{defin} 

\begin{thm}\label{equalsigns} 
$\wh \sigma(W_1(t), W_2(t))=\sigma(W_1(t), W_2(t))$. 
\end{thm}

\pp Given a vector 
$\e \in T$, we set $\e'=\e^{W_1}$ and $\e''=\e^{W_2}$. 
Choose a 
basis $\{\e_1, \ldots, \e_{n-1}\}$ of $T$ containing of 
the 
eigenvectors of $B$, i.e., $B\e_i=k_i\e_i$. Because of 
equality 
\eqref{oper} we have 
\begin{equation} 
\e''_i-\e'_i=(B\e_i)^S=(k_i\e_i)^S=k_i\e^S_i. 
\end{equation} 
We can and shall assume that the frame $\{\e_1, \ldots, 
\e_{n-1}\}$ gives the positive (local) orientation of $W_i 
(t), 
i=1,2$ at $x$. Take the polyvector 
$$ 
\mathfrak p: =\e'_1\wedge \cdots \wedge\e'_{n-1}\wedge 
{\mathbf v} 
\wedge \e''_1\wedge \cdots\wedge \e''_{n-1} 
$$ 
where ${\mathbf v}$ is the vector defined in $\ref{vector} 
$. Then 
$\mathfrak p \ne 0$ since the pair of propagations is 
assumed to 
be generic. Notice that 
$\mathfrak p$ gives us an orientation of $STM$, and we say 
that 
$\mathfrak p$ is positive if this orientation coincides 
with the 
original one, otherwise we say that $\mathfrak p$ is 
negative. 

According to Definition \ref{signcri}, the sign of 
$\mathfrak p$ is 
equal to 
$$ 
\eps(W_1(t), W_2(t))\sigma(W_1(t), W_2(t)). 
$$ 
So, we must prove that sign of the polyvector 
$\mathfrak p$ is equal to the sign of 
$$ 
(\det B)\sign(|{\mathbf v}_1|-|{\mathbf v}_2|). 
$$ 
To be definite, we assume that $|{\mathbf v}_1| > | 
{\mathbf v}_2|$ 
and prove that the sign of $\mathfrak p$ is equal to the 
sign of 
$\det B$. 

\m Since $\e''_i=\e'_i+k_i\e_i^S$ and $\det B=k_1\cdots k_ 
{n-1}$, 
we conclude that 
$$ 
\e'_1\wedge \cdots \wedge\e'_{n-1}\wedge {\mathbf v} 
\wedge 
\e''_1\wedge \cdots\wedge \e''_{n-1} 
=(\det B) \e'_1\wedge \cdots 
\wedge\e'_{n-1}\wedge {\mathbf v} \wedge \e^S_1\wedge 
\cdots\wedge 
\e^S_{n-1}. 
$$ 
So, it remains to prove that the polyvector 
$$ 
\e'_1\wedge \cdots 
\wedge\e'_{n-1}\wedge {\mathbf v} \wedge \e^S_1\wedge 
\cdots\wedge 
\e^S_{n-1} 
$$ 
is positive. 

Since $\e'_1\wedge \cdots \wedge\e'_{n-1}\wedge {\mathbf 
v} \wedge 
\e^S_1\wedge \cdots\wedge \e^S_{n-1}\ne 0$, we conclude 
that the 
family $\{\e^S_1,\ldots, \e^S_{n-1}\}$ generate $T_zS_x$. 
It is 
easy to see that the frame $\{\e^S_1, \ldots ,\e^S_{n-1}\} 
$ gives 
the original orientation of $S_x$, since the frame 
$\{\e_1, 
\ldots, \e_{n-1}\}$ gives the positive local orientation 
of each 
of the fronts at $x$. So, it remains to prove that the 
polyvector 
$$ 
(d\pr)(\e'_1\wedge \cdots \wedge\e'_{n-1}\wedge {\mathbf 
v})= 
\e_1\wedge \cdots \wedge\e_{n-1}\wedge (d\pr){\mathbf v} 
$$ 
is positive, i.e. that it gives the original orientation 
of $M$. 
Recall that the polyvector 
$$ 
\e_1\wedge \cdots \wedge\e_{n-1}\wedge ({\mathbf v}_1- 
{\mathbf v}_2) 
$$ 
is positive since $|{\mathbf v}_1|>|{\mathbf v}_2|$. 
Taking into 
account that the vector ${\mathbf v}_1-{\mathbf v}_2$ is 
orthogonal 
to $T$ and points into the direction of both propagations 
at $x$, it 
suffices to prove that $\langle {\mathbf v}_1-{\mathbf v} 
_2, 
(d\pr)({\mathbf v}\rangle > 0$. But this is clear because 
$W_1$ is faster 
then $W_2$ at $(x,t)$, and so the point $\pr(z_1)$ (in 
notation of \ref{vector}) 
is further then the point $\pr(z_2)$ from $T$. 
\qed 

\m 
{\em The following Proposition~$\ref{deformation}$ is 
useful, when calculating $\sign \det B$.\/} 
Let $\ov g $ be another Riemannian metric. Now we define 
$\ov H, \ov\pi, 
\ov\ell, \ov z$ and $\ov G$ as in \ref{setup} with respect 
to the metric $\ov 
g$. We also assume that $\ov\ell_1 
(U_1)\cap \ov\ell_2(U_2)=\ov z$. Note that in this case 
the image of $\ov 
\ell$ does not lie in $\Im \wt W$, but we can still define 
the operator $B=B(\ov g)$ as it is done before \propref 
{kerb}. 
We say that a metric $\ov g$ is {\it generic} if $d 
(\ov\ell_1)T_x\Im 
W_1\cap d(\ov\ell_2)T_x\Im W_2=0$. 

\begin{prop}\label{deformation} 
If the metric $\ov g$ is generic, then $\det B(\ov g)\ne 0 
$ and $\sign \det 
B(g)=\sign \det B(\ov g)$. 
\end{prop} 

\begin{rem} Probably the most useful case of this 
proposition is when $\ov g$ is chosen to be locally flat 
in the neighborhood of the dangerous tangency point and so 
that one of the fronts is (locally) a totally geodesic 
submanifold. 
\end{rem}

{\em Proof\/} of Proposition~\ref{deformation}. First, notice that the set of all Riemannian metrics 
is path connected 
because it is convex. 
Second, notice that $\det B(\wh g)\ne 0$ for every generic 
metric $\hat g$. (This 
is proved in the same way as \propref{kerb}.) 
Furthermore, the set of non-generic metrics has 
codimension $>1$ in the space of 
all  metrics. So, there exists a continuous family 
$g_t, t\in [0,1],$ of 
Riemannian metrics such $g_0=g, g_1=\ov g$ and each $g_t$ 
is generic. 
Clearly, $B(g_t)$ depends on $t$ continuously and $\det B 
(g_t)\ne 0$ since $g_t$ 
is generic. Thus the sign of $\det B(g_t)$ is the same for 
all $t$. 
\qed 

\begin{ex}[of calculation of $\sigma(\wt W_1(t), \wt W_2 
(t))$] Suppose that the fronts propagate as it is shown in 
Figure~\ref{definesign.fig}. 

\begin{figure}[htbp] 
\begin{center} 
\epsfxsize 12cm 
\hepsffile{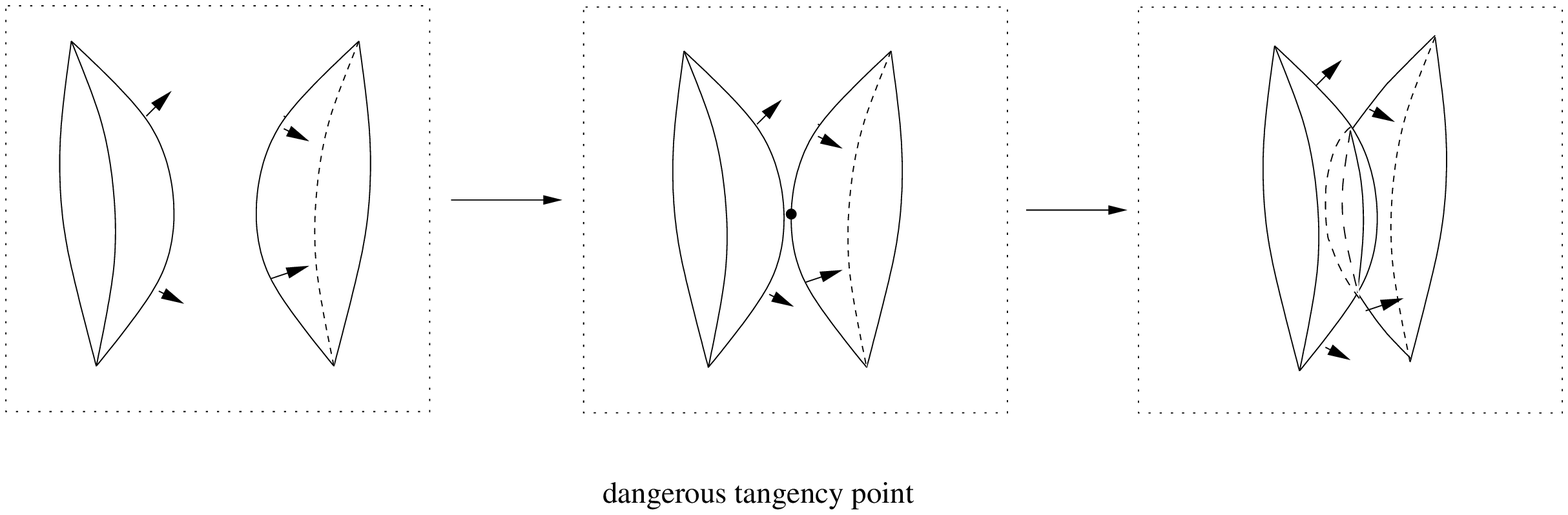} 
\end{center} 
\caption{}\label{definesign.fig} 
\end{figure} 
Let $W_1$ be the ``right'' front. Then, clearly, 
$$ 
\sigma(W_1(t),W_2(t))=-\eps(W_1(t),W_2(t)). 
$$ 
The negative sign appears because $W_2$ is faster 
the $W_1$. 
\end{ex} 

\section{Examples}\label{examples} 

To illustrate the usage of the affine linking invariant 
consider 
the following examples. 

\begin{ex} 
Here we show how to apply $\lk$ to determining the 
causality 
relation. Let $M$ be a smooth oriented $n$-dimensional 
manifold 
that is not an odd-dimensional homotopy sphere. Let $W_1, 
W_2$ be 
the wave fronts that originated on $M$ long time ago and 
were 
propagating according to the dangerous intersections free 
pair of 
propagations $\{P_1, P_2\}$.

Assume that the current picture of wave fronts $W_1(t), 
W_2 (t)$ 
is the one shown in Figure~\ref{example1causality.fig} 
with the 
velocity vectors normal to the two spheres shown in 
Figure~\ref{example1causality.fig}. (Note that after the 
contracting front 
contracts to a point it does not appear but rather everts 
and turns into an 
expanding spherical front.)

\begin{figure}[htbp] 
\begin{center} 
\epsfxsize 11cm 
\hepsffile{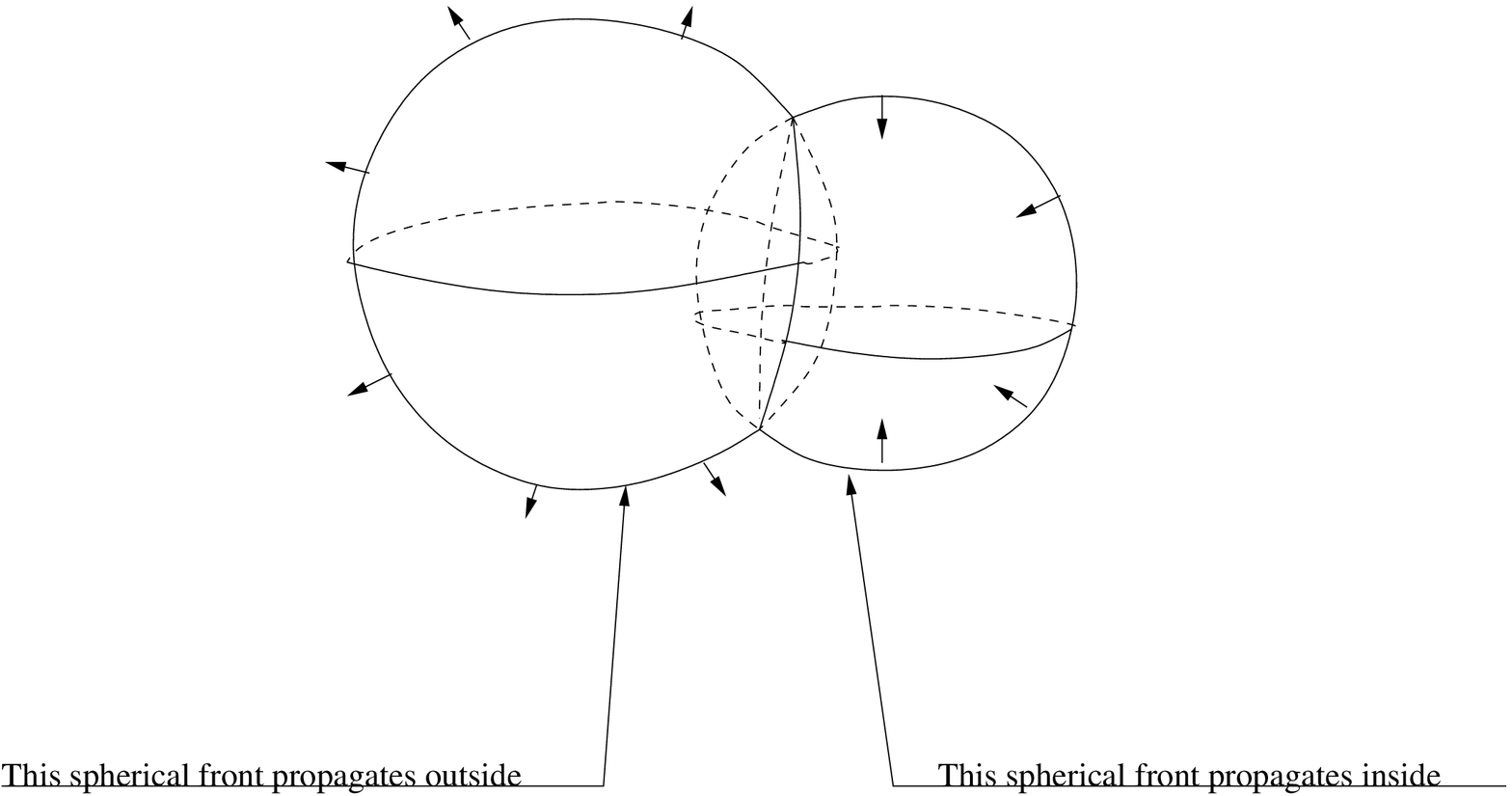} 
\end{center} 
\caption{}\label{example1causality.fig} 
\end{figure}

Then a straightforward calculation shows that 
$\lk(\wt W_1(t), \wt W_2(t))-\lk(V_1, V_2)=\pm 
1\neq 0$ (we used the notation as in Theorem~\ref{main}), 
and thus the first wave front reached the birth point of 
the second front before the second front originated. 
(The sign of $\pm 1$ in 
this example depends on which of the two fronts shown in 
Figure~\ref{example1causality.fig} is $W_1$ in the case 
where $n$ is odd and is 
always a plus sign when $n$ is even.) 

This seems to demonstrate that $\lk$ is a very powerful 
invariant 
because in this case we know neither the propagation laws 
nor when 
and where the fronts originated. In fact, in this example 
we can 
make this conclusion even without the knowledge of the 
topology of 
$M$ outside of the depicted part of it. 
\end{ex}

\begin{ex} 
Here we show how to apply $\lk$ to estimating of 
the number of times the wave front passed 
though a given point between the two moments of time. 

Assume that we have a wave front $W$ that propagates on 
$M$ and 
that $M$ is not an odd-dimensional homotopy sphere. 

Assume that at a certain moment of time the picture of the 
wave front was the one shown in Figure~\ref 
{example2causality.fig}.a and later 
it developed into the shape shown in Figure~\ref 
{example2causality.fig}.b. 
(The Figure~\ref{example2causality.fig}.b depicts a sphere 
that can be obtained 
from the trivially embedded sphere 
by passing three times through a 
point and by creation of some singularities far away from 
$x$.) 

\begin{figure}[htbp] 
\begin{center} 
\epsfxsize 10cm 
\hepsffile{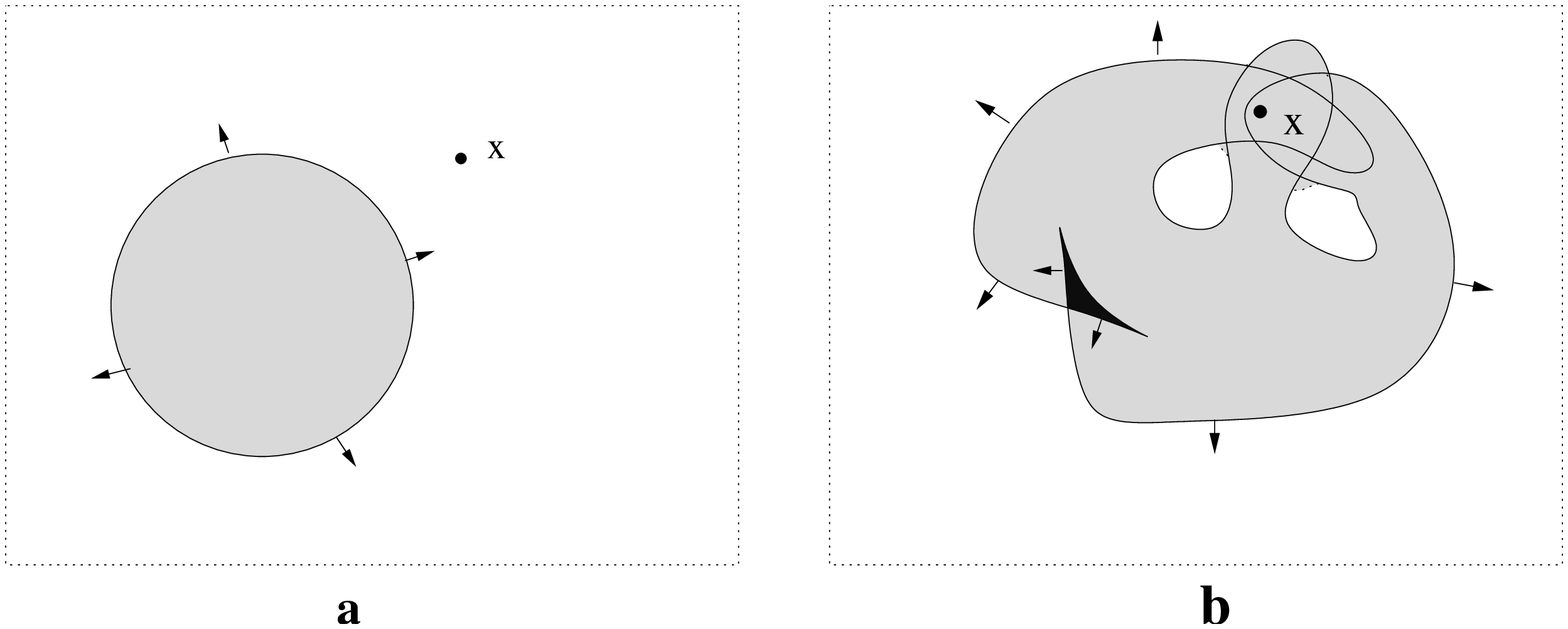} 
\end{center} 
\caption{}\label{example2causality.fig} 
\end{figure}

Let $\wt W(t_1), \wt W(t_2): S^{n-1} \rightarrow STM$ be 
the 
liftings of the fronts 
shown in Figure~\ref{example2causality.fig}.a and~b 
respectively. A straightforward 
calculation shows that 
$\lk(\wt W(t_2), \eps_x)-\lk(\wt W(t_1), \eps_x)=3\in\A(M) 
$ 
for every map $\eps_x:S^{n-1} \to STM$ as in \eqref{eps}. 

Thus if the dimension of the ambient manifold is even, or 
if $\pi_1(M)$ is 
infinite, then (see~\ref{usefulprop}) $W$ did pass at 
least three times through 
the point $x$ 
between the time moments shown in Figure~\ref 
{example2causality.fig}.a and~\ref{example2causality.fig} 
b. 
Once again, this conclusion does not depend on the 
topology 
of $M$ outside of 
the part of it depicted in Figure~\ref 
{example2causality.fig}, on the time 
passed between the two pictures taken, and on the 
propagation law. 
\end{ex}

{\bf Acknowledgments:} The first author was supported by 
the 
free-term research money from the Dartmouth College. The 
second author 
was supported by the free-term research money from the 
University 
of Florida, Gainesville and by MCyT, 
projects BFM 2002-00788 and MCyT BFM2003-02068/MATE, 
Spain. 

The authors are very thankful to Robert Caldwell, Jose Natario, and 
Jacobo 
Pejsachowicz for useful discussions.

\end{document}